\documentclass{article} % For LaTeX2e
\usepackage{times}
\usepackage{fullpage}
\usepackage{natbib}
\usepackage{amsmath,amsfonts,bm}

\def\eqref#1{equation~\ref{#1}}
\def\1{\bm{1}}

\DeclareMathAlphabet{\mathsfit}{\encodingdefault}{\sfdefault}{m}{sl}
\SetMathAlphabet{\mathsfit}{bold}{\encodingdefault}{\sfdefault}{bx}{n}

\newcommand{\R}{\mathbb{R}}

\usepackage{hyperref}
\usepackage{url}

\usepackage{ifthen,algorithm,algorithmic}
\usepackage{amssymb,amsfonts,amsmath,latexsym,dsfont}
\usepackage{tikz,pgfplots,pgflibraryplotmarks}
\usetikzlibrary{calc}
\usepackage{graphicx}
\usepackage{mathrsfs}
\usepackage{algorithm,algorithmic}
\usepackage{hyperref}
\usepackage{multirow}

\graphicspath{{images/}{./}}

\usepackage{boxedminipage}

\usepackage{pgf,pgfarrows} %% pdf-drawing package
\usepackage{subfigure} %% pdf-drawing package

\newcommand{\FrameboxA}[2][]{#2}
\newcommand{\Framebox}[1][]{\FrameboxA}

\newcommand{\bfC}{{\bf C}}

\newcommand{\bfF}{{\bf F}}

\newcommand{\bfI}{{\bf I}}

\newcommand{\bfK}{{\bf K}}

\newcommand{\bfM}{{\bf M}}

\newcommand{\bfW}{{\bf W}}

\newcommand{\bfc}{{\bf c}}
\newcommand{\bfe}{{\bf e}}

\newcommand{\bfx}{ {\bf x}}
\newcommand{\bfy}{ {\bf y}}

\newdimen\iwidth
\newdimen\iheight

\newcommand{\CO}{{\cal O}}

\title{Low-Cost Parameterizations of Deep Convolutional Neural Networks}

\author{Eran Treister\thanks{Department of Computer Science, Ben-Gurion University of the Negev, Beer Sheva, Israel. ({\tt erant@cs.bgu.ac.il, [sharmic,sapirza]@post.bgu.ac.il})} \and Lars Ruthotto\thanks{Department of Mathematics and Computer Science, Emory University, Atlanta, Georgia, USA. ({\tt lruthotto@emory.edu})}$^{\;}$\thanks{Xtract Inc., Vancouver, Canada}
\and Michal Sharoni\footnotemark[1] \and Sapir Zafrani\footnotemark[1] \and Eldad Haber\thanks{ Earth and Ocean Sciences Department, UBC, Vancouver, Canada. ({\tt ehaber@eos.ubc.ca})}$^{\;}$\footnotemark[3]\;}

\begin{document}

\maketitle
\begin{abstract}
Convolutional Neural Networks (CNNs) filter the input data using a series of spatial convolution operators with compactly supported stencils and point-wise nonlinearities.
Commonly, the convolution operators couple features from all channels.
For wide networks, this leads to immense computational cost in the training of and prediction with CNNs.
In this paper, we present novel ways to parameterize the convolution more efficiently, aiming to decrease the number of parameters in CNNs and their computational complexity.
We propose new architectures that use a sparser coupling between the channels and thereby reduce both the number of trainable weights and the computational cost of the CNN.
Our architectures arise as new types of residual neural network (ResNet) that can be seen as discretizations of a Partial Differential Equations (PDEs) and thus have predictable theoretical properties. Our first architecture involves a convolution operator with a special sparsity structure, and is applicable to a large class of CNNs. Next, we present an architecture that can be seen as a discretization of a diffusion reaction PDE, and use it with three different convolution operators. We outline in our experiments that the proposed architectures,  although considerably reducing the number of trainable weights, yield comparable accuracy to existing CNNs that are fully coupled in the channel dimension.
\end{abstract}

\section{Introduction}

Convolutional Neural Networks (CNNs)~\citep{LeCun1990} are among the most effective machine learning approaches for processing structured high-dimensional input data and are indispensable in,  e.g., in recognition tasks involving speech~\citep{RainaEtAl2009} and image data~\citep{KrizhevskySutskeverHinton2012}.
The essential idea behind CNNs is to replace some or all of the affine linear transformations in a neural network by convolution operators that are typically parameterized using small-dimensional stencils.
This has a number of benefits including the increase of computational efficiency of the network due to the sparse connections between features, and a considerable reduction in the number of weights since stencils are shared across the whole feature map~\citep{Goodfellow:2016wc}.

In practical applications of CNNs, the features can be grouped into channels whose number is associated with the width of the network.
This gives one several opportunities to define interactions between the different channels.
Perhaps, the most common approach in CNNs is to fully couple features across channels~\citep{Gu:2018id,Goodfellow:2016wc,KrizhevskySutskeverHinton2012}.
Following this approach, the number of convolution operators at a layer is proportional to the product of the number of input and output channels.
Given that performing convolutions is often the computationally most  expensive part in training of and prediction with CNNs and the number of channels is large in many applications, this scaling can be problematic for wide architectures or high-dimensional data.
Another disadvantage of this type of architecture is the number of weights. Indeed, for deep neural networks, the number of weights that are associated with a wide network can easily reach millions and beyond. This makes the deployment of such networks challenging, especially on devices with limited memory.

In this paper, we propose four novel ways to parameterize CNNs more efficiently, based on ideas from Partial Differential Equations (PDEs).
Our goal is to dramatically reduce the number of weights in the networks and the computational costs of training and evaluating the CNNs. One ides, similar to \cite{howard2017mobilenets}, is to use spatial convolutions for each channel individually and add $1\times 1$ convolutions to impose coupling between them. Our architectures are motivated by the interpretation of residual neural networks (ResNets)~\citep{he2016deep,he2016identity} as time-dependent nonlinear PDEs~\citep{ruthotto2018deep}. More specifically, we consider a simple Reaction-Diffusion (RD) model, that can model highly nonlinear processes. We derive new architectures by discretizing this continuous model, using $1\times 1$ convolutions as a reaction term, together with cheap forms of a spatial convolution, that are similar to a depth-wise convolution in the number of parameters and cost. This spatial convolution acts similarly to a linear diffusion term that smooths the feature channels. Since the networks we propose originate in continuous models they have distinct theoretical properties that can be  predicted using the standard theory of ODEs and PDEs \citep{AscherPetzoldODEs}.

Our first approach is designed to be employed in any existing CNN layer with equal number of input and output channels. We simply replace the traditional fully coupled convolution with a \emph{linear} sum of depth-wise and $1\times1$ convolution, like a mask that can be placed on a traditional convolution in any existing CNN. Our second ``explicit'' RD architecture applies the operators separately with a non-linear activation function operating only following the $1\times1$ convolution, as the non-linear reaction part of the diffusion reaction equation. The third architecture is more unique.
To improve the stability of the forward propagation and increase the spatial coupling of the features, we propose a semi-implicit scheme for the forward propagation through the network. Unlike traditional CNN operators, the semi-implicit scheme applies an inverse of the depth-wise (block diagonal) convolution preceded by a non-linear step involving the $1 \times 1$ convolution. This way, the scheme couples {\em all} the pixels in the image in one layer, even though we are using a depth-wise $3\times 3$ convolution. The inverse of the convolution operator can be efficiently computed using Fast Fourier Transforms (FFT) and over the channels and kernels.

The last idea is to replace the depth-wise convolution structure with a circulant connectivity between the channels.
This is motivated by the interpretation of the features as tensors and follows the definition of an efficient tensor product in~\citep{Kernfeld:2015fo} whose associated tensor singular value decomposition has been successfully used for image classification in~\citep{NewmanEtAl2017}.
The circulant structure renders the number of trainable convolution stencils proportional to the width of the layer.
Using periodic boundary conditions in the other feature dimensions, this convolution can be computed efficiently by extending the FFT-based convolutions in~\citep{Mathieu:2013wa,Vasilache:2014wh} along the channel dimension, which reduces the cost from ${\cal O}(c^2)$ to  ${\cal O}(c \log c)$ where $c$ is the number of channels.

Table~\ref{tab:costs} compares the number of weights and the computational complexity associated with the forward propagation through a layer for the standard and reduced architectures. In the table we assume that the explicit RD architecture is directly computed without using FFT, but the FFT-based implementation, which is necessary for the implicit scheme, can also be used for the explicit one.

Our architectures pursue a similar goal than the recently proposed MobileNet architectures that are also based on a mix of $1\times1$ and ``depth-wise'' convolutions \citep{howard2017mobilenets,sandler2018mobilenetv2}.
% This depth-wise convolution is a block diagonal convolution operator that separately operates on each feature channel without any connection across the channels.
The MobileNet architecture involves with significantly less parameters, and in particular avoids the fully coupled convolutions, except for $1\times1$ convolutions which are cheaper in both computational cost and number of parameters.
What sets our work apart from these architectures is that our architectures can be seen as discretization of PDEs, which allows to control their stability and offers new ways for their analysis.

\begin{table}[t]
	\caption{Cost comparison of different convolution layers for an image with $n$ pixels, stencil of size $m\times m$, and $c$ input and output channels. RD denotes a reaction-diffusion architecture.} \label{tab:costs}
    \centering
	\begin{tabular}{rrr}
    \hline
		               & no. of weights & computational costs \\
					  % \cmidrule(r){2-2} \cmidrule(r){3-3}
		fully-coupled  & $\CO(m^2\cdot c^2)$ & $\CO(n \cdot m^2 \cdot c^2)$ \\
		RD explicit  & $\CO(m^2\cdot c + c^2)$ & $\CO(n (m^2\cdot c + c^2))$ \\
		RD implicit  & $\CO(m^2\cdot c + c^2)$ & $\CO(n(\log(n)\cdot c + c^2))$ \\
        RD circulant      & $\CO(m^2\cdot c + c^2)$ & $\CO((n c) \log (cn))$ \\
    \hline
	\end{tabular}
\end{table}

The remainder of the paper is organized as follows.
We first describe the mathematical formulation of the supervised classification problem with deep residual neural networks used in this paper.
Subsequently, we propose the novel parameterizations of CNNs, describe their efficient implementation, and their computational costs.
We perform experiments using the CIFAR10, CIFAR 100, and STL10 datasets and demonstrate that the performance of the new architectures, despite a considerable reduction in the number of trainable weights, is comparable to residual neural networks using fully-coupled convolutions.
Finally, we summarize and conclude the paper.

\section{Mathematical Formulation} % (fold)
\label{sec:math}

In this section, we introduce our main notation and briefly describe the mathematical formulation of the supervised classification problem used in this paper, which is based on~\citep{Goodfellow:2016wc}.
For brevity we restrict the discussions to images although techniques derived here can also be used for other structured data types such as speech or video data.

Given a set of training data consisting of image-label pairs, $\{(\bfy^{(k)}, \bfc^{(k)})\}_{k=1}^s \subset \R^{n_f}\times \R^{n_c}$ and a residual neural network (ResNets)~\citep{he2016deep,he2016identity}, our goal is to find network parameters $\theta \in \R^n$ and weights of a linear classifier defined by $\bfW, \mu \in\R^{n_c}$ such that
\begin{equation*}
	\bfc^{(k)} \approx S(\bfW \bfy^{(k)}(\theta) +  \mu), \quad \text{ for all } \quad k=1,2,\ldots,s,
\end{equation*}
where $S$ is a softmax hypothesis function and $\bfy^{(k)}(\theta)$ denotes the output features of the neural network applied to the $k$th image.
As common, we model the learning problem as an optimization problem aiming at minimizing a regularized empirical loss function
\begin{equation*}
	\min_{\theta,\bfW,\mu} \frac1s \sum_{k=1}^s L(S(\bfW \bfy^{(k)}(\theta) + \mu), \bfc^{(k)}) + R(\theta,\bfW,\mu),
\end{equation*}
where in this paper $L$ is the cross entropy and $R$ is a smooth and convex regularization function.
Solving the optimization problem is not the main focus of the paper, and we employ standard variants of the stochastic gradient descent algorithm (see the original work of~\cite{RobbinsMonro1951} and the survey of~\citep{bottou2016optimization}).

As common in other deep learning approaches, the performance of the image classification hinges upon designing an effective the neural network, i.e., the relation between the input feature $\bfy^{(k)}$ and its filtered version $\bfy^{(k)}(\theta)$.
The goal in training is to find a $\theta$ that transforms the features in a way that simplifies the classification problem.
In this paper, we restrict the discussion to convolutional ResNets, which have been very successful for image classification tasks.
As pointed out by recent works~\citep{HaberRuthotto2017,Chang2017Reversible,HaberHolthamRuthotto2017,E2017,ChaudhariEtAl2017,ruthotto2018deep} there is a connection between ResNets and partial differential equations (PDEs), which allows one to analyze the stability of a convolution ResNet architecture.

We consider a standard ResNet as a baseline architecture. For a generic example $\bfy_0$, the $j$th ``time step''  of the network reads
\begin{equation}\label{eq:Resnet}
    \bfy_{j+1} = \bfy_{j} + \bfF(\theta_{j},\bfy_j), \quad \text{ for all } \quad j=0,1,\ldots,N-1,
\end{equation}
where $\theta_{j}$ are the weights associated with the $j$th step. This step can be seen as a discretization of the initial value problem
\begin{equation}\label{eq:PDE}
    \partial_{t} \bfy(t) = \bfF(\theta(t),\bfy(t)), \quad \bfy(0) = \bfy_0, \quad t \in [0,T].
\end{equation}
In a simple ResNet, the nonlinear term in \eqref{eq:Resnet} and~\eqref{eq:PDE} usually reads
\begin{equation}\label{eq:F}
\bfF(\theta,\bfy) = \bfK(\theta^{(2)})^\top \sigma( {\cal N}(\bfK(\theta^{(1)}) \bfy)),
\end{equation}
where $\sigma(x) = \max\{x,0\}$ denotes a element-wise rectified linear unit (ReLU) activation function, the weight vector is partitioned into $\theta^{(1)}$ and $\theta^{(2)}$ that parameterizes the two linear operators $\bfK(\theta^{(1)})$ and $\bfK(\theta^{(1)})$, and $\mathcal{N}$ denotes a normalization layer that may have parameters as well (omitted here).

In CNNs the linear operator $\bfK$ in~\eqref{eq:F} is defined by combining spatial convolution operators, which gives a rich set of modeling options by, e.g., specifying boundary conditions, padding, strides, and dilation~\citep{Goodfellow:2016wc}.
Our focus in this work is the coupling between different feature channels.
Assuming that the input feature $\bfy$ can be grouped into $c$ channels, the most common choice is to use full coupling in across the channels.
As an example, let $\bfy$ consist of four channels and the number of output channels be four as well. Then, the operator $\bfK(\theta)$ is a four by four block matrix consisting of convolution matrices $\bfC$ parametrized by the different stencils that comprise $\theta$
\begin{equation}\label{eq:CNNfc}
	\bfK(\theta)
	=
	\left(
		\begin{array}{llll}
			\bfC\left(\theta^{(1)}\right) &  \bfC\left(\theta^{(2)}\right) &  \bfC\left(\theta^{(3)}\right)&  \bfC\left(\theta^{(4)}\right)\\
			\bfC\left(\theta^{(5)}\right) &  \bfC\left(\theta^{(6)}\right) &  \bfC\left(\theta^{(7)}\right)&  \bfC\left(\theta^{(8)}\right)\\
			\bfC\left(\theta^{(9)}\right) &  \bfC\left(\theta^{(10)}\right) & \bfC\left(\theta^{(11)}\right)& \bfC\left(\theta^{(12)}\right)\\
			\bfC\left(\theta^{(13)}\right) & \bfC\left(\theta^{(14)}\right) & \bfC\left(\theta^{(15)}\right)& \bfC\left(\theta^{(16)}\right)\\		
		\end{array}
	\right).
\end{equation}
Hence, applying a convolution operator  going from $c_{\rm in}$ to $c_{\rm out}$ channels requires $\mathcal{O}( c_{\rm in} \cdot c_{\rm out})$ convolutions.
Since practical implementations of convolutional ResNets often use hundreds of channels, this coupling pattern leads to large computational costs and to millions of parameters.
Hence, we define more efficient coupling strategies in the next section.
% section related_work (end)

\section{Low-Cost Parameterizations of Convolution Operators} % (fold)
\label{sec:low_cost_parameterizations}
In this section, we present novel ways to parameterize the convolution operators in CNNs more efficiently.
The first convolution operator is a simple sum of a depth-wise and $1\times1$ convolution, which can be thought of as a masked version of \eqref{eq:CNNfc}.
The next two networks are discretizations of a new type of ResNet that are inspired by reaction-diffusion PDEs, where depth-wise spatial convolution operator is used as a diffusion process.
The last approach imposes a block circulant structure on the diffusion operator in the reaction-diffusion PDE.
We also provide detailed instructions on how to implement the proposed architectures efficiently.

\subsection{The depth-wise and $1\times 1$ convolutions}
Most of our new architectures are defined using two types of operators. One is the depth-wise (block diagonal) operator which operates on each channel separately. Omitting the step index $j$, the operator in matrix form is given by
\begin{equation}\label{eq:Kbd}
	\bfK_{\rm dw}(\theta) = 	\left(
		\begin{array}{llll}
			\bfC\left(\theta^{(1)}\right) &   &  &  \\
			  & \bfC\left(\theta^{(2)}\right) &&  \\
			 &  & \bfC\left(\theta^{(3)}\right)& \\
			 & & & \bfC\left(\theta^{(4)}\right)\\		
		 \end{array}
	\right).
\end{equation}

Another building block is the fully-coupled $1\times1$ convolution operator that couples different channels but introduces no spatial filtering. For $\theta \in \R^{16}$, we denote such an operator as
\begin{equation}\label{eq:1X1}	
	\bfM(\theta) = \left(
		\begin{array}{llll}
			\theta_{1}  & \theta_{5}    &\theta_{9}    &\theta_{13}   \\
			\theta_{2}  & \theta_{6}    &\theta_{10}   &\theta_{14}   \\
			\theta_{3} & \theta_{7}    &\theta_{11}   &\theta_{15}   \\
			\theta_{4}  & \theta_{8}    &\theta_{12}   &\theta_{16}   \\
		\end{array}
	\right) \otimes \bfI,
\end{equation}
where $\otimes$ denotes the Kronecker product, and $\bfI$ is an identity matrix.
Note that $\bfM(\theta)$ models a reaction between the features in different channels at a given pixel, but introduces no spatial coupling. In MobileNets \citep{howard2017mobilenets}, the operators \eqref{eq:Kbd} and \eqref{eq:1X1} are used interchangeably in separate neural network layers, which also include non-linear activation and batch-normalization for each layer separately. We note that the depth-wise convolution is related to the ``2D-filter'' structurally sparse convolutions in \citep{wen2016learning}. There, however, the authors allow the kernels to be all over the matrix and not only on the diagonal, and choose them via group $\ell_1$ penalty.

Our first idea to simplify the coupled convolution in \eqref{eq:CNNfc}, is to combine the diagonal and off-diagonal weights into one operator that is multiplied as a standard matrix:
\begin{eqnarray}\label{eq:Mask}	
	\bfK_{\rm lm}(\theta^{(1)},\theta^{(2)}) = \bfK_{\rm dw}(\theta^{(1)})  + \bfM(\theta^{(2)}).
\end{eqnarray}
This is a masked version of the original convolution \eqref{eq:CNNfc}, with the mask leaving just the block diagonal and $1\times1$ convolution terms (note that the diagonal part of $\bfM(\theta^{(2)})$ can be ignored as the same coefficients appear also in $\bfK_{\rm dw}(\theta^{(1)})$).
This type of convolution can be used instead of the standard operators in CNNs, as long as the number of input and output channels are equal.

% paragraph circulant_kernel (end)
\subsection{Networks based on the Reaction-Diffusion equation}

%As an alternative to the second-order hyperbolic forward propagation in~\eqref{eq:wave} we introduce a new class of CNNs based on the reaction-diffusion like equation
As an alternative to the architecture in equations~\ref{eq:Resnet}-\ref{eq:F} we introduce a new class of CNNs based on the reaction-diffusion like equation
\begin{eqnarray}\label{eq:rdPDE}
	\partial_t \bfy(t) =  - \bfK_{\rm dw}\left(\theta^{(1)}\right)^\top \bfK_{\rm dw}\left(\theta^{(1)}\right) \bfy(t)
                       + \sigma\left({\cal N}\left(\bfM\left(\theta^{(2)}\right) \bfy(t)\right)\right), \qquad \bfy(0) = \bfy_0.
\end{eqnarray}
Such equations have been used to model highly nonlinear processes such as pattern formation in complex chemical, biological and social systems. They typically lead to interesting patterns which suggests that they are highly expressive. Therefore, using such models for learning is very intriguing and a natural extension for standard ResNets.

\subsubsection{Explicit Reaction Diffusion CNN.}

In the simplest, ``explicit'', discretization of the RD equation we use a ResNet structure of
\begin{eqnarray}\label{eq:explicit}
\bfy_{j+1} = \bfy_j + h  \left(-\bfK_{\rm dw}\left(\theta_j^{(1)}\right)^\top \bfK_{\rm dw}\left(\theta_j^{(1)}\right)\bfy_j + \sigma\left({\cal N}(\bfM(\theta_j^{(2)}) \bfy_j)\right) \right),
\end{eqnarray}
%with a special choice for each of the matrices $\bfK_1$ and $\bfK_2$.
%$$\bfK_2(\theta^{(2)}) = \bfI - \bfK_{\rm dw}\left(\theta^{(2)}\right)^\top \bfK_{\rm dw}\left(\theta^{(2)}\right),$$
%an identity minus a symmetric positive definite block diagonal matrix, and
%$$\bfK_1(\theta^{(1)}) = \bfM(\theta^{(1)}),$$
where the time step $h>0$ is chosen sufficiently small.
The first symmetric and positive (semi) definite convolution operates as a diffusion on each channel separately, while the second term---the $1\times 1$ convolution--- models a reaction between the features in different channels at a given pixel, without introducing spatial coupling.

Both types of linear operators in~\eqref{eq:explicit} can be implemented efficiently.
The fully coupled part in $\bfM$ is identical to the standard $1 \times 1$ convolution, and can be computed by a single call to a matrix-matrix multiplication BLAS routine ({\tt gemm}) without any need to manipulate or copy the data.
For 2D convolutions with filter size $m$ the cost of this operation is a factor $m^2$ cheaper.
The depth-wise operator $\bfK_{\rm dw}$  can be computed directly by applying standard convolutions on each of the channels. However, it can also be computed using FFT.
Similar as above, each operator $\bfC$ in \eqref{eq:Kbd} can be evaluated by
\begin{equation}\label{eq:explicit_comp}
\tilde\bfC \bfy = \bfF_2^{-1} \left((\bfF_2 (\tilde\bfC\bfe_1))  \odot  (\bfF_2 \bfy) \right),
\end{equation}
where  $\bfF_2$ and $\bfF_2^{-1}$ are the 2D FFT and inverse FFT, respectively, and $\bfe_1$ is the first standard basis vector. The cost of this computation scales linearly with the number of channels compared with the number of channels square of the standard fully connected convolution. The convolution is applied using the batched 2D FFT  routines in the library {\tt cufft}.

\subsubsection{Implicit Reaction Diffusion CNN. } % (fold)
\label{par:diagional_kernel}
Our second type of CNN may be seen as an ``implicit'' version of the previous convolutional layer, which is known to be a stable way to discretize the forward propagation~\eqref{eq:rdPDE}. Here, we use a semi-implicit time-stepping
\begin{eqnarray}\label{eq:implicit}
	\bfy_{j+1} = \left(\bfI + h \bfK_{\rm dw}(\theta_j^{(1)})^\top\bfK_{\rm dw}(\theta_j^{(1)})\right)^{-1}
                \left(\bfy_j + h \sigma ({\cal N}(\bfM(\theta_j^{(2)}) \bfy_j)) \right).
\end{eqnarray}
{The better stability of the implicit equation stems from the \emph{unconditionally} bounded spectral radius of the first operator, specifically
$$
\rho\left(\left(\bfI + h \bfK_{\rm dw}(\theta)^\top\bfK_{\rm dw}(\theta)\right)^{-1}\right) < 1, \quad \forall \theta.
$$
Since this matrix is part of the Jacobian of the forward step (with respect to the input), the stability properties of the implicit forward propagation \citep{ruthotto2018deep} are better than those of its explicit counterpart. This behaviour is well known in the context of time-dependent PDEs.}
%In addition to improving the stability of the forward propagation,
In addition, the implicit step has another advantage. It yields a global coupling between pixels in only one application, which allows features in one side of the image to impact features in other side.
The coupling decays away from the center of the convolution and is strongly related to the Green's Function of the associated convolution.

We exploit the special structure of $\bfK_{\rm dw}$ to efficiently solve the linear system in \eqref{eq:implicit} using ${\cal O}(c\cdot n \log(n))$ operations.
While, in general, inverting a $K\times K$ matrix inversion requires a ${\cal O}(K^3)$ operation, the depth-wise kernel $\bfK_{\rm dw}$ is block diagonal and therefore, the inverse is computed by solving  each of the $c$ blocks individually and in parallel. Since each diagonal block is a convolution, its inverse can be computed efficiently using a 2D FFT when assuming periodic boundary conditions.  To this end we use the formula
$$ (h \bfC^{\top}\bfC + \bfI)^{-1} \bfy = \bfF_2^{-1} \left(
(h |\bfF_2 \bfC\bfe_1|^2 + 1)^{-1} \odot  (\bfF_2 \bfy) \right). $$
The FFT operations are essentially identical to the ones in \eqref{eq:explicit_comp}, and hence have similar cost.

%An example for such a coupling for random weights  is plotted in Figure~\ref{figInvConv}.
% \begin{figure}
% \begin{center}
% \includegraphics[width=3cm]{invConv.jpg}
% \caption{The inverse convolution operator of a random $3\times 3$ kernel on a $32\times 32$ grid. The weights are non-zero all over the grid, but decay away from the center of the convolution.\label{figInvConv}}
% \end{center}
% \end{figure}
%As we see next, this model can be also evaluated efficiently compared to the standard convolution model.

\subsubsection{Block Circulant Convolution as the diffusion.} % (fold)
\label{par:circulant_kernel}
Another way to increase the computational efficiency of CNNs is based on the interpretation of the image data as a 3D tensor whose third dimension represents the channels. Following the notion of the tensor product in~\citep{Kernfeld:2015fo}, we define using the block circulant operator
\begin{equation}\label{eq:CNNcir}
	\bfK_{\rm circ}(\theta)
	=
	\left(
		\begin{array}{llll}
			\bfC\left(\theta^{(1)}\right) & \bfC\left(\theta^{(2)}\right) & \bfC\left(\theta^{(3)}\right) & \bfC\left(\theta^{(4)}\right)\\
			\bfC\left(\theta^{(4)}\right) & \bfC\left(\theta^{(1)}\right) & \bfC\left(\theta^{(2)}\right) & \bfC\left(\theta^{(3)}\right)\\
			\bfC\left(\theta^{(3)}\right) & \bfC\left(\theta^{(4)}\right) & \bfC\left(\theta^{(1)}\right) & \bfC\left(\theta^{(2)}\right)\\
			\bfC\left(\theta^{(2)}\right) & \bfC\left(\theta^{(3)}\right) & \bfC\left(\theta^{(4)}\right) & \bfC\left(\theta^{(1)}\right)\\
		\end{array}
	\right).
\end{equation}
Using the associated tensor SVD has shown promising results~\citep{NewmanEtAl2017} on the MNIST data set.
% \tred{[ET: Is the next sentence supposed to be here?]}Results have been extended to CNNs using the single-layer Hamiltoninan network described in~\cite{HaberRuthotto2017}, which is less expressive than the ResNet model we consider here.
% As we see next, the computational models above can be computed in a significantly less operations compared with the standard fully coupled convolution.
We assume periodic boundary conditions (potentially requiring padding).
Under this assumption, a matrix-vector product between the block circulant matrix with circulant blocks, $\bfK_{\rm circ}$ and a feature vector,
can be done using a Fast Fourier Transform (FFT), i.e.,
$$\bfK_{\rm circ} \bfy = \bfF_3^{-1} ( (\bfF_3 (\bfK_{\rm circ} \bfe_1)) \odot (\bfF_3 \bfy)). $$
Here, $\bfF_3$ is  a 3D Fast Fourier Transform (FFT),  $\bfF_3^{-1}$ is its inverse (see, e.g.,~\citep{HansenNagyOLeary2006} for details).
The computational complexity of this product in proportional to $(nc) \log (nc)$ where
$n$ is the number of pixels, compared to the order of $m^2 n c^2$ of the fully coupled convolution. It also requires much less parameters. While standard convolution requires $m^2 \cdot c^2$ variables our convolution requires only $m^2 \cdot c$ variables that we save in a 3D array.

\section{Experiments} % (fold)
\label{sec:numerical_examples}

We experimentally compare the architectures proposed in this paper to the ResNet and MobileNet architectures using the CIFAR-10, CIFAR100~\citep{krizhevsky2009learning} and STL-10~\citep{coates2011analysis} data sets.
Our primary focus is on showing that similar accuracy can be achieved using a considerably smaller number of weights in the networks. All our experiments are performed with the PyTorch software \citep{paszke2017automatic}.
%Since our focus is on the computational aspects of each of these networks we do not employ techniques such as data augmentation or excessive tuning of hyperparameters and employ the same basic architecture for all experiments, comparing the fully coupled models with our new models.

\paragraph{Base Architecture. }
We use the same base architecture in all our numerical experiments, which is a slightly modified version of the one described in~\citep{Chang2017Reversible}. We use three network sizes, but the overall structure is identical between them. Our goal is to use simple, standard, and rather small networks, and show that the new architectures can lead to a performance comparable to a standard ResNet architecture using even less parameters.

Our networks consist of several blocks, that are preceded by an opening layer.
This opening layer is a convolutional layer with a $5\times5$ convolution  that increases the number of channels from 3 to 32 or 48, depending on the network. This is followed by a batch normalization, and a ReLu activation. Then, there are several blocks (three or four), each consisting of a ResNet based part with four steps that varies between the different experiments except the ReLu activation and batch normalization. The architectures for the series of steps are:
\begin{itemize}
\item ResNet - a step with two fully coupled convolutions as defined in \eqref{eq:Resnet}-\eqref{eq:F}.
\item MobileNet - a two layer neural network similar to \citep{howard2017mobilenets}
\begin{equation}
\hat\bfy = \sigma({\cal N}(\bfK_{\rm bd}(\theta^{(1)})\bfy_{j}));\quad \bfy_{j+1} = \sigma({\cal N}(\bfM(\theta^{(2)})\hat\bfy)).
\end{equation}
\item LinearMix - a ResNet step using \eqref{eq:Mask} as operators.
\item Explicit / implicit RD - the architectures in \eqref{eq:explicit}, \eqref{eq:implicit} respectively.
\item Circular RD - the architectures in \eqref{eq:explicit} with $\bfK_{\rm circ}$ in \eqref{eq:CNNcir} instead of $\bfK_{\rm dw}$.
\end{itemize}

\begin{table}[t]
	\caption{Classification results }
    \centering
	\begin{tabular}{|l|cc|cc|cc|}
\hline
                  & \multicolumn{2}{c|}{CIFAR10} & \multicolumn{2}{c|}{CIFAR100}& \multicolumn{2}{c|}{STL10}\\
Architecture      & Network   & val. acc.  & Network &  val. acc.  & Network & val. acc. \\
\hline
 ResNet           & A  (1.5M)     &  93.1\%  &  B (3.4M) & 71.7\% & A(1.5M)& 74.9\%\\
 ResNet	          &  B (3.5M)     &  93.0\%  &  C (6.3M) & 69.9\%  &B(3.5M)& 75.3\% \\
 MobileNet        &  A  (101K)      & 89.5\%  &  B (251K) & 65.6\%& A(101K)& 74.9\%\\
 MobileNet        &  B  (216K)     & 91.6\%   &  C (423K) & 61.9\%& B(216K)& 77.2\%\\
 LinearMix    &  A  (195K)      &  91.3\% & B (456K)& 67.9\%& A(195K)& 75.6\%\\
 LinearMix    &  B  (422K)     &  92.1\%  & C (789K)& 69.2\% &B (422K)&75.6\% \\
 Exp. RD            & A   (101K)      &  88.9\%  & B (250K)& 66.1\% &A (101K)&74.6\% \\
 Exp. RD            & B   (216K)      &  90.6\%  & C (423K)& 65.2\% &B (216K)&75.9\% \\
 Imp. RD            & A  (101K)      &  88.7\%   & B (250K)& 64.6\% &A(101K)& 73.8\%\\
 Imp. RD            & B  (216K)      &  90.3\%   & C (423K)& 64.9\%& B(216K)&73.4\% \\
 Circ. RD           & A  (101K)      &  86.0\%  & B (250K)& 60.2\%& A(101K)& 69.6\%\\
 Circ. RD           & B  (216K)      &  88.0\%  & C (423K)& 60.0\% &B (216K)&70.5\% \\
\hline
	\end{tabular}\label{tab:Classification}
\end{table}

Each series of steps is followed by a single ``connecting'' layer that takes the images and concatenate them the same images multiplied with a depth-wise convolution operator and batch-normalization:
$$
\bfx \leftarrow \cal N([\bfx ; \bfK_{\rm dw}(\theta)\bfx]).
$$
This doubles the number of channels, and following this we have an average pooling layer with down-sample the images by a factor of 2 at each dimension. We have also experimented with other connecting layers, such as a more standard $1\times1$ convolution either followed by pooling or with strides, leading to similar results. We use three networks that differ in the number of channels:
$$
A: 32-64-128 \quad B: 48-96-192 \quad C:32-64-128-256
$$
The last block consists of a pooling layer that averages the image intensities of each channel to a single pixel and we use a fully-connected linear classifier with softmax and cross entropy loss.

As the number of parameters is typically small, we do not use regularization for training the networks. For training the networks we use the ADAM optimizer \citep{kingma2014adam}, with its default parameters. We run 300 epochs and reduce the learning rate by a factor or 0.5 every 60 epochs, starting from 0.01. We used a mini-batch size of 100. We also used standard data augmentation, i.e., random resizing, cropping and horizontal flipping.

\paragraph{Data Sets: }
The CIFAR-10 and CIFAR100 datasets~\citep{krizhevsky2009learning} consists of 60,000 natural images of size $32\times32$ with labels assigning each image into one of ten categories (for CIFAR10) or 100 categories (for CIFAR100).
The data are split into 50,000 training and 10,000 test images. The STL-10 dataset~\citep{coates2011analysis} contains $13,000$ color-images each of size $96\times96$ that are divided into $5,000$ training and $8,000$ test images that are split into the ten categories.

Our classification results are given in Table \ref{tab:Classification}. The results show that our different architectures are in par and in some cases better than other networks. The theoretical properties of our architectures can be explained by the standard theory of ODEs and PDEs which makes them more predictable given small perturbations in the network parameters (for example, truncation errors) or noise in the data \citep{ruthotto2018deep}. The architecture used is computationally efficient and the efficiency increases as the number of channels increases. For many problems where the number of channels is in the thousands, our approach can yield significant benefits compared with other architectures.

\subsection{Computational Performance}

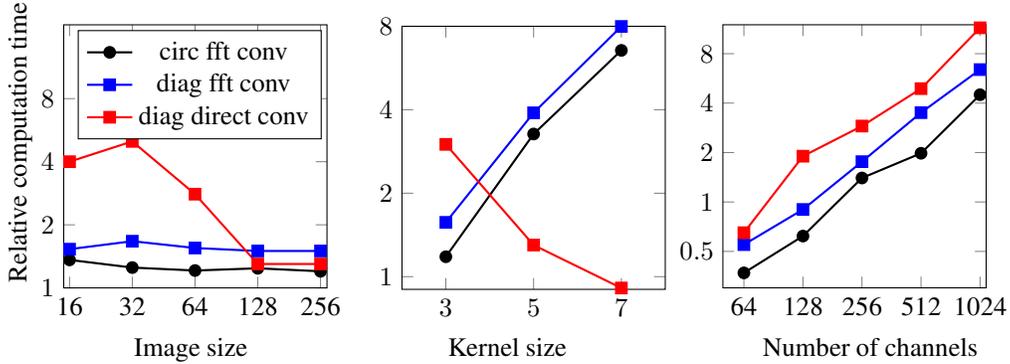
\begin{figure*}[t]
	\begin{center}
		\newcommand{\rottext}[1]{\rotatebox{90}{\hbox to 30mm{\hss  #1\hss}}}
		\iwidth=35mm
		\iheight=35mm
		\begin{tabular}{@{}c@{}c@{}c@{}c@{}}
		%&  $f=0.5$ & $f_i=\{0.5,0.75\}$ \\
		\rottext{\hspace{50pt} Relative computation time}   & % \begin{tikzpicture}
%     \begin{axis}[
%     ymode=log,
%     xmode=log,
%     log ticks with fixed point
%     ]
% \addplot table {
% 9   1.75439
% 14  0.71429
% 18  0.41667
% 19  0.33333
% 34  0.11765
% };
% \end{axis}
% \end{tikzpicture}

\begin{tikzpicture}
\begin{axis}[%
width=\iwidth,
height=\iheight,
at={(0,0)},
scale only axis,
log ticks with fixed point,
xmin=15,
xmax=275,
xtick={16,32,64,128,256},
ytick={1,2,4,8},
ymin=1.0,
ymax=18,
ymode=log,
xmode=log,
xticklabels from table={Im1.dat}{input},
%yminorticks=true
% legend style={at={(0.03,0.03)},anchor=south west,legend cell align=left,align=left,draw=white!15!black}
]
\addplot [color=black,thick,solid,mark=*]
  table[row sep=crcr]{%
16   1.36  \\
32   1.25  \\
64   1.21 \\
128   1.24 \\
256   1.2 \\
};
\addplot [color=blue,thick,solid,mark=square*]
table[row sep=crcr]{%
16   1.53  \\
32   1.67  \\
64   1.55 \\
128   1.5 \\
256   1.5 \\
};
\addplot [color=red,thick,solid,mark=square*]
table[row sep=crcr]{%
16   4.0  \\
32   5.0  \\
64   2.8 \\
128   1.3 \\
256   1.3 \\
};

\addlegendentry{circ fft conv}
\addlegendentry{diag fft conv}
\addlegendentry{diag direct conv}

\end{axis}
\end{tikzpicture}%
 &  \hspace{8pt}\begin{tikzpicture}
\begin{axis}[%
width=\iwidth,
height=\iheight,
at={(0,0)},
scale only axis,
xmin=2,
xmax=8,
ymin=0.9,
xtick={3,5,7},
ytick={1,2,4,8},
ymax=8,
ymode=log,
log ticks with fixed point,
%xmode=log
% xlabel={iterations},
%yminorticks=true
% legend style={at={(0.03,0.03)},anchor=south west,legend cell align=left,align=left,draw=white!15!black}
]
\addplot [color=black,thick,solid,mark=*]
  table[row sep=crcr]{%
3   1.18  \\
5   3.27  \\
7   6.55 \\
};
\addplot [color=blue,thick,solid,mark=square*]
  table[row sep=crcr]{%
3   1.57  \\
5   3.9  \\
7   7.99 \\
};
\addplot [color=red,thick,solid,mark=square*]
  table[row sep=crcr]{%
3   3.0  \\
5   1.3  \\
7   0.91 \\
};

\end{axis}
\end{tikzpicture}%
& \hspace{2pt}% \begin{tikzpicture}
%     \begin{axis}[
%     ymode=log,
%     xmode=log,
%     log ticks with fixed point
%     ]
% \addplot table {
% 9   1.75439
% 14  0.71429
% 18  0.41667
% 19  0.33333
% 34  0.11765
% };
% \end{axis}
% \end{tikzpicture}

\begin{tikzpicture}
\begin{axis}[%
width=\iwidth,
height=\iheight,
at={(0,0)},
scale only axis,
log ticks with fixed point,
xmin=50,
xmax=1100,
xtick={64,128,256,512,1024},
ytick={0.25,0.5,1,2,4,8},
ymin=0.3,
ymax=12,
ymode=log,
xmode=log,
xticklabels from table={Ch1.dat}{input},
%yminorticks=true
% legend style={at={(0.03,0.03)},anchor=south west,legend cell align=left,align=left,draw=white!15!black}
]
\addplot [color=black,thick,solid,mark=*]
  table[row sep=crcr]{%
64   0.37  \\
128   0.62  \\
256   1.40 \\
512   1.98 \\
1024   4.5 \\
};
\addplot [color=blue,thick,solid,mark=square*]
table[row sep=crcr]{%
 64   0.55  \\
 128   0.90  \\
 256   1.76 \\
 512   3.5 \\
 1024   6.4 \\
};
\addplot [color=red,thick,solid,mark=square*]
table[row sep=crcr]{%
 64   0.65  \\
 128   1.9  \\
 256   2.9 \\
 512   4.9 \\
 1024   11.5 \\
};

\end{axis}
\end{tikzpicture}%
\\
		& Image size  &  Kernel size &\hspace{20pt} Number of channels\\
		\end{tabular}
	\end{center}
	\caption{Runtime ratio between {\tt cudnn}'s fully coupled convolution, depth-wise convolution and our implementation of the depth-wise (explicit or implicit) and circular convolutions. Except the measured parameter, the default parameters for all tests are: batch size: 64, image size: $64^2$, number of channels: 256, kernel size $3 \times 3$.}
	\label{fig:CompTime}
\end{figure*}

In this section we compare the runtime of the forward step of our FFT-based convolutions relatively to the runtime of the fully coupled convolution in \eqref{eq:CNNfc}, which is computed using the {\tt cudnn} package version 7.1. This package is used in all of the GPU implementations of CNN frameworks known to us. Both the circular and depth-wise convolutions are implemented using {\tt cufft} as noted above. For the direct depth-wise convolution, we use PyTorch's implementation using {\tt groups}.

The experiments are run on a Titan Xp GPU. We report the relative computation time
$$
\frac{time(\mbox{Fully coupled using {\tt cudnn}})}{time(\mbox{Our implementation})},
$$
so we wish that the ratio is large. We note that the convolutions \eqref{eq:CNNcir} and \eqref{eq:Kbd} can also be applied using {\tt cudnn} by manually forming the convolution kernel (zero-filled, or circularly replicated), and hence we can use {\tt cudnn} in cases where the time ratio is smaller than 1.

Our tests are reported in Fig.~\ref{fig:CompTime}. The presented runtime ratio was calculated based on the total time of 100 convolutions. The left most graph presents the execution time ratio with respect to the image size. As expected, the execution time ratio does not depend on the image size in all convolutions, except for the direct method at small scales, which may be faster due to efficient memory access. The middle graph presents the execution time ratio with respect to the stencil size. Here, since both of our implementations apply the FFT on a zero-padded kernel weights, their execution time is independent of the kernel size, and the time ratio compared to the direct {\tt gemm}-based cudnn implementation improves as the kernel grows. The FFT-based implementation is favorable if one wishes to enrich the depth-wise convolution models with wider stencils. The left-most graph presents the execution ratio with respect to the number of channels. Here we can clearly see that the execution ratio is linear in the number of channel, with a ratio of 1 achieved at about 200 channels for the FFT-based implementations. The direct convolution has a better constant, but the overall complexity is similar. Clearly, the considered convolutions are more favorable for wide networks.

\section{Discussion and Conclusion} \label{sec:discussion}
We present four new convolution models with the common goal of reducing the number of parameters and computational costs of  CNNs.
To this end, we propose alternative ways to the traditional full coupling of channels, and thereby obtain architectures that involve fewer expensive convolutions, avoid redundancies in the network parametrization, and thereby can be deployed more widely.
Our work is similar to that of \citep{howard2017mobilenets,sandler2018mobilenetv2}. However, our unique angle is the close relation of our architectures to continuous models given in terms of PDEs that are well understood. This highlights stability of our CNNs and paves the way toward more extensive theory.

Our numerical experiments for image classification show that the new architectures can be almost as effective as more expensive fully coupled CNN architectures. We expect that our architectures will be able to replace the traditional convolutions in classification of audio and video, and also in other tasks that are treated with CNNs. It is important to realize that our new architectures become even more advantageous for 3D or 4D problems, e.g., when analyzing time series of medical or geophysical images. In these cases, the cost of each convolution  is much more expensive and the computational complexity makes using 3D CNNs difficult. Here, also the number of weights imposes challenges when using computational hardware with moderate memory.

% \newpage
% section discussion_and_conclusion (end)
%\bibliographystyle{iclr2019_conference}
\bibliographystyle{plain}

\bibliography{LowCostBib}

\end{document}